\documentclass[12pt]{article}
\usepackage{amssymb, amsmath, amsthm, mathrsfs}
\usepackage{hyperref}

\def\ddt{\frac{\d}{\d t}}

\def\d{{\rm d}}

\def\e{{\rm{e}}}

\def\R{\mathbb R}
\def\lap{\Delta}
\def\grad{\nabla}
\def\div{\nabla\cdot}

\def\Re{\R}
\def\half{{\textstyle\frac{1}{2}}}

\def\qas{\quad{\rm as}\quad}

\def\qqfa{\qquad\mbox{for all}\qquad}
\def\qqand{\qquad\mbox{and}\qquad}
\def\qand{\quad\mbox{and}\quad}
\def\qqwith{\qquad\mbox{with}\qquad}

\def\be#1{\begin{equation}\label{#1}}
\def\ee{\end{equation}}
\def\bea{\begin{eqnarray*}}
\def\eea{\end{eqnarray*}}

\def\({\left(}
\def\){\right)}
\def\[{\left\lbrack}
\def\]{\right\rbrack}

\newtheorem{theorem}{Theorem}

\newtheorem{lemma}[theorem]{Lemma}

\newtheorem{corollary}[theorem]{Corollary}

\date{September, 18, 2006}

\begin{document}

\title{{\it A posteriori} regularity of the three-dimensional Navier-Stokes
equations from numerical computations{\thanks {{\bf To Appear in:}
{\it Journal of Mathematical Physics.}}}}

\author{Sergei I. Chernyshenko,\\
{\small Aeronautics and Astronautics, School of Engineering Sciences,}\\
{\small University of Southampton,}\\
{\small Highfield, Southampton, SO17 1BJ. UK.}\\ \\
Peter Constantin,\\
{\small Department of Mathematics, University of Chicago,}\\
{\small 5734 University Avenue, Chicago, IL 60637. USA.}\\ \\
James C. Robinson,\\
{\small Mathematics Institute, University of Warwick, Coventry, CV4
7AL. UK.}\\ \\
Edriss S. Titi\\
{\small Department of Mathematics}\\
{\small and Department of Mechanical \& Aerospace Engineering,}\\
{\small University of California, Irvine, CA 92697-3875. USA.}\\
{\small Also: Department of Computer Science \& Applied Mathematics,}\\
{\small Weizmann Institute of Science, Rehovot 76100. Israel.}}

\maketitle

\newpage

\section*{Abstract}

%{Navier-Stokes equations, rigorous computation, Galerkin method}
In this paper we consider the r\^ole that numerical computations --
in particular Galerkin approximations -- can play in problems
modelled by the 3d Navier-Stokes equations, for which no rigorous
proof of the existence of unique solutions is currently available.
We prove a robustness theorem for strong solutions, from which we
derive an {\it a posteriori} check that can be applied to a
numerical solution to guarantee the existence of a strong solution
of the corresponding exact problem.

We then consider Galerkin approximations, and show that {\it if} a
strong solution exists the Galerkin approximations will converge
to it; thus if one is prepared to assume that the Navier-Stokes
equations are regular one can justify this particular numerical
method rigorously.

Combining these two results we show that if a strong solution of
the exact problem exists then this can be verified numerically
using an algorithm that can be guaranteed to terminate in a finite
time.

We thus introduce the possibility of rigorous computations of the
solutions of the 3d Navier-Stokes equations (despite the lack of
rigorous existence and uniqueness results), and demonstrate that
numerical investigation can be used to rule out the occurrence of
possible singularities in particular examples.
%\end{abstract}

\newpage

\section{Introduction}

The Navier-Stokes equations are the fundamental model of fluid flow.
However, no proof of the long-time existence of unique solutions of
these equations is currently available. The importance of this
long-standing problem has recently been highlighted by its inclusion
as one of the Clay Foundation's Millennium Problems.

Despite this mathematical obstacle, the equations are used as the
basis of numerical calculations both in theoretical investigations
of turbulence and for industrial applications. This paper addresses
the validity of such numerical computations given the lack of
appropriate rigorous existence and uniqueness results.

One might expect that without the existence of a unique solution
there is no hope of guaranteeing that a numerical approximation is
really an `approximation' in any meaningful sense, since it is not
clear what is being approximated. However, here we prove three
results that enable numerical and exact solutions to be related in a
rigorous way.

First we provide an explicit check that can be applied to a
numerical solution to guarantee that it is approximating a solution
of the exact problem. In particular, this check implies the
existence of a unique solution of the exact problem over the same
time interval as the calculation. Thus the validity of a numerical
solution can be verified rigorously by a simple {\it a posteriori}
condition.

We then turn to the particular example of Galerkin approximations.
We show that given the {\it assumption} of the existence of a
sufficiently smooth unique solution of the exact problem, the
solutions obtained by the Galerkin method will converge to this
exact solution: If one is prepared to take on trust that the
Navier-Stokes equations are a meaningful model of fluid flow,
numerical experiments (at least those using the Galerkin method) can
be justified mathematically.

Finally, these two results can be combined to show that the
existence of a sufficiently smooth solution can be verified
numerically using an algorithm that can be guaranteed to terminate
within a finite time (provided that the solution exists).

Although we concentrate here on spatial discretization (mainly via
the spectral Galerkin method), it is relatively straightforward to
consider fully discrete methods, and these are also addressed where
relevant.

In order to discuss our results in a little more detail, we
introduce here the model and its abstract setting. We will consider
the 3d Navier-Stokes equations for an incompressible fluid
\be{realnse}
\frac{\partial u}{\partial t}-\nu\lap u+(u\cdot\grad)u+\grad
p=\tilde f(t)\qqwith\div u=0
\ee
on a periodic domain $Q=[0,L]^3$ with the additional (convenient)
restriction that $\int_Q u=\int_Q \tilde f=0$. However, we would
expect similar results to hold for more physically realistic
boundary conditions (e.g.~Dirichlet boundary conditions on a bounded
domain).

In order to recast the equations in their functional form (for full
details see the monographs by Temam (1977), Constantin \& Foias
(1988), or Robinson (2001)) we let $\mathscr H$ be the collection of
all divergence-free smooth periodic vector-valued functions on $Q$
with zero average, and set
\bea
H&=&\mbox{closure of }{\mathscr H}\mbox{ in }[L^2(Q)]^3\\
V&=&\mbox{closure of }{\mathscr H}\mbox{ in }[H^1(Q)]^3.
\eea
We use $|\cdot|$ and $(\cdot,\cdot)$ for, respectively, the norm and
inner product in $[L^2(Q)]^3$.

Denoting by $\Pi$ the orthogonal projection of $[L^2(Q)]^3$ onto
$H$, we apply $\Pi$ to (\ref{realnse}) and obtain
\be{nse_one}
\frac{\d u}{\d t}+\nu Au+B(u,u)=f(t),
\ee
where $A$ is the Stokes operator $Au=-\Pi\lap u$ (in fact $Au=-\lap
u$ when $u\in D(A)$ for the periodic case which we consider here),
$B(u,u)=\Pi[(u\cdot\grad)u]$, and $f=\Pi\tilde f$. The pressure term
disappears since gradients are orthogonal (in $L^2$) to
divergence-free fields.

We note that in the space periodic case we have $D(A^{m/2})=H^m\cap
V$; for simplicity we denote this space by $V^m$, and equip it with
the natural norm
$$
\|u\|_m=|A^{m/2}u|
$$
which is equivalent to the norm in the Sobolev space $H^m$.

In section 2 we present a very simple ODE lemma that forms the basis
of what follows. Section 3 recalls some classical estimates on the
nonlinear term $B(u,u)$ in the spaces $V^m$ ($m\ge 3$), and
reproduces the proof of a classical regularity result.

Section 4.1 proves the robustness of strong solutions: We show in
theorem \ref{robust} for $m\ge 3$ that if the initial data $u_0\in
V^m$ and forcing $f(t)\in L^2(0,T;V^{m-1})\cap L^1(0,T;V^m)$ give
rise to a strong solution, then so do `nearby' initial data $v_0\in
V^m$ and forcing functions $g\in L^2(0,T;V^{m-1})$, under the
explicit condition that
\begin{eqnarray*}
&&\|u_0-v_0\|_m+\int_0^T\|f(s)-g(s)\|_m\,\d
s\\
&&\qquad\qquad<\frac{1}{c_mT}\exp\[-c_m\int_0^T\(\|u(s)\|_m
+\|u(s)\|_{m+1}\)\,\d s\]
\end{eqnarray*}
(note that this depends only on the solution $u$). A similar result
holds for the Euler equations under slightly modified hypotheses,
and this is shown in section 4.2.

As corollary \ref{theretheyare} in section 5 we use the observation
that any sufficiently smooth numerically computed solution $u$ is
the exact solution of the Navier-Stokes equation for some
appropriate forcing function to turn the previous robustness result
into an {\it a posteriori} test for the existence of a strong
solution based on numerical observations. Leaving the precise
definition of `sufficiently smooth' to the formal statement of the
result, we show that if $u$ is a `good' approximation to the
solution of
\be{nseintro}
\frac{\d v}{\d t}+\nu Av+B(v,v)=f(t)\qqwith v(0)=v_0\in V^m
\ee
in the sense that
\begin{eqnarray*}
&&\|u(0)-v_0\|_m+\int_0^T\left\|\frac{\d u}{\d t}(s)+\nu
Au(s)+B(u(s),u(s))-f(s)\right\|_m\,\d s\\
&&\qquad\qquad<\frac{1}{c_mT}\exp\[-c_m\int_0^T\(\|u(s)\|_m+\|u(s)\|_{m+1}\)\,\d
s\]
\end{eqnarray*}
then $v(t)$ must be a strong solution of (\ref{nseintro}) on $[0,T]$
with $v\in L^\infty(0,T;V^m)\cap L^2(0,T;V^{m+1})$. Crucially, this
condition depends only on the numerical solution $u$ (and the given
initial data and forcing, $v_0$ and $f(t)$).

In section 6 we consider the convergence of approximations to
solutions of the Navier-Stokes obtained via a Fourier Galerkin
method. Here we {\it assume} the existence of a sufficiently regular
exact solution (in fact we take $u_0\in V^m$ and $f\in L^2(0,T;V^m)$
and assume the existence of a strong solution), and demonstrate in
theorem \ref{galerkin_thm} that under this condition, its Galerkin
approximations $u_n$ do indeed converge to the correct limit in both
$L^\infty(0,T;V^m)$ and $L^2(0,T;V^{m+1})$.

Finally we combine our {\it a posteriori} test and the convergence
of Galerkin approximations in theorem \ref{verify} to show that if,
for a given choice of initial data $u_0\in V^m$ and forcing $f\in
L^1(0,T;V^m)\cap L^2(0,T;V^{m-1})$ (where $m\ge3$), a strong
solution does exist for some time interval $[0,T]$, then this can be
verified computationally in a finite number of steps. (Essentially
we show that for a sufficiently large Galerkin calculation we can
ensure that the {\it a posteriori} test of corollary
\ref{theretheyare} must be satisfied.)

It should be noted that we do not aim here to prove results that are
optimal with regard to the regularity of solutions, but rather to
take sufficient regularity (initial conditions in the Sobolev space
$H^m$ with $m\ge 3$) that the arguments are at their most
straightforward. Similar results for initial conditions in $H^1$ and
$H^2$ are possible, and will be presented elsewhere.

\section{A simple ODE lemma}

The following simple lemma, after that in Constantin (1986), is
central to all the results that follow. Although all solutions of
the equation
$$
\dot y=\delta(t)+\alpha y^2\qqwith y(0)=y_0>0\qand\alpha>0
$$
blow up in a finite time, if $y_0$ and $\delta(t)$ are sufficiently
small (in appropriate senses) then the solution can be guaranteed to
exist on $[0,T]$.

\begin{lemma}\label{lemma}
Let $T>0$ and $\alpha>0$ be constants, and let $\delta(t)$ be a
non-negative continuous function on $[0,T]$. Suppose that $y$
satisfies the differential inequality
\be{ineq}
\frac{\d y}{\d t}\le \delta(t)+\alpha y^2\qqwith y(0)=y_0\ge 0
\ee
and define
$$
\eta=y_0+\int_0^T\delta(s)\,\d s.
$$
Then
\be{est1}
y(t)\le\frac{\eta}{1-\alpha\eta t}
\ee
while $0\le t\le T$ and the right-hand side is finite. In particular
if
\be{condition}
\alpha\eta T<1
\ee
then $y(t)$ remains bounded on $[0,T]$, and clearly
$y(t)\rightarrow0$ uniformly on $[0,T]$ as $\eta\rightarrow0$.
\end{lemma}

\begin{proof}
Observe that $y(t)$ is bounded by $Y(t)$ where
$$
\frac{\d Y}{\d t}=\alpha Y^2\qqwith Y(0)=y_0+\int_0^T\delta(s)\,\d
s,
$$
which yields (\ref{est1}) and the remainder of the lemma follows
immediately.\end{proof}

\section{A classical regularity result}

In this section we prove a slight variant of a classical regularity
result which can be found as theorem 10.6 in Constantin \& Foias,
1988.

First, we recall the following bounds on the nonlinear term, which
we will use repeatedly (see Constantin \& Foias (1988) for a proof):
\begin{align}
\|B(u,v)\|_m&\le c_m\|u\|_m\|v\|_{m+1}&m&\ge 2,\label{Buv}\\
|(B(w,v),A^mw)|&\le c_m\|v\|_{m+1}\|w\|_m^2&m&\ge 2,\label{mm1}\\
|(B(v,w),A^mw)|&\le c_m\|v\|_m\|w\|_m^2&m&\ge 3.\label{mm}
\end{align}

\begin{theorem}\label{reggie}
Let $u\in L^\infty(0,T;V)\cap L^2(0,T;V^2)$ be a strong solution
of the 3d Navier-Stokes equations with $u_0\in V^m$ and $f\in
L^2(0,T;V^{m-1})$. Then in fact
$$
u\in L^\infty(0,T;V^m)\cap L^2(0,T;V^{m+1})
$$
and
\be{dudt}
\d u/\d t\in L^2(0,T;V^{m-1}).
\ee
\end{theorem}

\begin{proof}
We give a formal argument which can be made rigorous using the
Galerkin procedure.  The proof is inductive, supposing initially
that $u\in L^2(0,T;V^k)$ for some $k\le m$. Taking the inner
product of the equation with $A^ku$ we obtain
$$
\frac{1}{2}\frac{\d}{\d t}\|u\|_k^2+\nu\|u\|_{k+1}^2\le
|(B(u,u),A^ku)|+(f,A^ku)
$$
and so using (\ref{mm1}) (valid here for $k\ge 2$)
$$
\frac{1}{2}\frac{\d}{\d t}\|u\|_k^2+\nu\|u\|_{k+1}^2\le
c_k\|u\|_k^2\|u\|_{k+1}+\|f\|_{k-1}\|u\|_{k+1}.
$$
Therefore
\be{backtome}
\frac{\d}{\d
t}\|u\|_k^2+\nu\|u\|_{k+1}^2\le\frac{c_k^2}{\nu}\|u\|_k^4+\frac{\|f\|_{k-1}^2}{\nu},
\ee
Dropping the term $\nu\|u\|_{k+1}^2$ we have
$$
\frac{\d}{\d
t}\|u\|_k^2\le\left(\frac{c_k^2\|u\|_k^2}{\nu}\right)\|u\|_k^2+\frac{\|f\|_{k-1}^2}{\nu};
$$
It now follows from the Gronwall inequality that our assumption
$u\in L^2(0,T;V^k)$ implies that $u\in L^\infty(0,T;V^k)$.

Returning to (\ref{backtome}) and integrating between $0$ and $T$
we obtain
\be{intbound}
\nu\int_0^T\|u(s)\|_{k+1}^2\,\d
s\le\|u(0)\|_k^2+\frac{c_k^2}{\nu}\int_0^T\|u(s)\|_k^4\,\d
s+\frac{1}{\nu}\int_0^T\|f(s)\|_{k-1}^2\,\d s,
\ee
which shows in turn that $u\in L^2(0,T;V^{k+1})$.

Since by assumption $u\in L^2(0,T;V^2)$, the first use of the
induction requires $k=2$, for which inequality (\ref{mm1}) is
valid: we can therefore conclude by induction that $u\in
L^\infty(0,T;V^m)\cap L^2(0,T;V^{m+1})$.
Finally, since
$$
\frac{\d u}{\d t}=-\nu Au-B(u,u)+f
$$
and each term on the right-hand side is contained in
$L^2(0,T;V^{m-1})$, the bound on the time derivative in
(\ref{dudt}) follows.\end{proof}

%Note that (\ref{intbound}) shows in particular that boundedness in
%$L^\infty(0,T;V^m)$ of a sequence of solutions $u_n$ of problems
%with $f_n$ bounded in $L^2(0,T;V^{m-1})$ implies that $u_n$ is
%also bounded in $L^2(0,T;V^{m+1})$.

We note here that it follows from this theorem that if $u_0$ and $f$ are smooth (in $V^m$ for all $m$) then so is the solution. Different techniques (due to Foias \& Temam, 1989) can be used to show that the solution is analytic in the space variable (i.e.~in a certain Gevrey class) provided that the data is.

\section{Robustness of strong solutions}

\subsection{The Navier-Stokes equations}

Using lemma \ref{lemma} we show that if the 3d Navier-Stokes
equations have a sufficiently smooth strong solution for given
initial data $u_0$ and forcing $f$ then they also have a strong
solution for close enough data. The argument is based closely on
that in Constantin (1986) which, given the existence of a strong
solution of the Euler equations, deduces the existence of strong
solutions for the Navier-Stokes equations for small enough
$\nu>0$. (The same argument is used in Chapter 11 of the monograph by Constantin \& Foias (1988)).

We choose to state our primary result for sufficiently smooth strong
solutions, namely those corresponding to initial data in $V^m$ with
$m\ge 3$. This enables us to use all the inequalities
(\ref{Buv}--\ref{mm}) and thereby obtain a relatively simple
`closeness' condition in (\ref{condition2}). A similar approach
works with strong solutions that have the minimal required
regularity ($u\in L^\infty(0,T;V)\cap L^2(0,T;V^2)$) but the results
are less elegant; these results will be presented in a future paper.

Note that while here we concentrate on the robustness of solutions
defined on finite time intervals, a result valid for all $t\ge0$
given a solution $u$ for which $\int_0^\infty\|u(s)\|_1^4\,\d
s<\infty$ (i.e.~which decays appropriately as $t\rightarrow\infty$)
has been obtained by Ponce et al.~(1993) for the particular case
$f=0$.

\begin{theorem}\label{robust}
Let $m\ge 3$ and let $u$ be a strong solution of the 3d
Navier-Stokes equations
$$
\frac{\d u}{\d t}+\nu Au+B(u,u)=f(t)\qqwith u(0)=u_0\in V^m
$$
and $f\in L^2(0,T;V^{m-1})\cap L^1(0,T;V^m)$. Then if $g\in
L^2(0,T;V^{m-1})$ and
\begin{eqnarray}
&&\|u_0-v_0\|_m+\int_0^T\|f(s)-g(s)\|_m\,\d
s\nonumber\\
&&\qquad\qquad<\frac{1}{c_mT}\exp\[-c_m\int_0^T\(\|u(s)\|_m
+\|u(s)\|_{m+1}\)\,\d s\]\label{condition2}
\end{eqnarray}
the solution $v$ of
$$
\frac{\d v}{\d t}+\nu Av+B(v,v)=g(t)\qqwith v(0)=v_0\in V^m
$$
is a strong solution on $[0,T]$ and is as regular as $u$.
\end{theorem}

We remark here that throughout this paper we consider only
Leray-Hopf weak solutions, i.e.~weak solutions satisfying the energy
inequality starting from almost every time.

\begin{proof}
Standard existence results guarantee that $v$ is a strong solution
on some time interval $[0,T^*)$. If $T^*$ is maximal then
$\|u(t)\|_1\rightarrow\infty$ as $t\rightarrow T^*$; clearly we also
have $\|u(t)\|_m\rightarrow\infty$ as $t\rightarrow T^*$. We suppose
that $T^*\le T$ and deduce a contradiction.

While $v$ remains strong, the assumption that $v_0\in V^m$ and
that $g\in L^2(0,T;V^{m-1})$ allows one to use the regularity
results of theorem \ref{reggie} to deduce that
$$
v\in L^\infty(0,T';V^m)\cap L^2(0,T';V^{m+1})
$$
for any $T'<T^*$; we also have $\d v/\d t\in L^2(0,T';V^{m-1})$
(cf.~the argument in the proof of theorem \ref{reggie}). It also
follows from theorem \ref{reggie} that the solution $u$ enjoys
similar regularity on $[0,T]$.

The difference $w=u-v$ satisfies
\be{star}
\frac{\d w}{\d t}+\nu Aw+B(u,w)+B(w,u)+B(w,w)=f-g
\ee
with $w(0)=w_0=u_0-v_0$. On $[0,T^*)$ we know that $w$ is
sufficiently regular that
\be{langlerangle}
\langle A^{m/2}\frac{\d w}{\d
t},A^{m/2}w\rangle=\frac{1}{2}\,\frac{\d}{\d t}|A^{m/2}w|^2.
\ee
We can therefore take the inner product with $A^mw$ and obtain,
using (\ref{mm1}) and (\ref{mm}),
\be{tooyoung}
\frac{1}{2}\ddt{}\|w\|_m^2\le
c_m\|u\|_m\|w\|_m^2+c_m\|u\|_{m+1}\|w\|_m^2+c_m\|w\|_m^{3}+\|f-g\|_m\|w\|_m.
\ee
Dividing by $\|w\|_m$ yields\footnote{Of course, one has to worry
here whether $\|w(t)\|_m$ is zero. However, if $\|w(t_0)\|_m=0$ for
some $t_0$ then $w(t_0)=0$, and then the uniqueness of strong
solutions in class of Leray-Hopf weak solutions implies that
$u(t)=v(t)$ for all $t\ge t_0$, i.e.~$w(t)=0$ for all $t\in[t_0,T]$,
a contradiction.}
\bea
\ddt{}\|w\|_m&\le&\|f-g\|_m+c_m\|u\|_m\|w\|_m+c_m\|u\|_{m+1}\|w\|_m+c_m\|w\|_m^2\\
&\le&\|f-g\|_m+c_m(\|u\|_m+\|u\|_{m+1})\|w\|_m+\|w\|_m^2.
\eea

We multiply by $\exp(-c_m\int_0^t\|u(s)\|_m+\|u(s)\|_{m+1}\,\d s)$
and consider
$$
y(t)=\|w(t)\|_m\exp\[-c_m\int_0^t\(\|u(s)\|_m+\|u(s)\|_{m+1}\)\,\d
s\],
$$
for which we obtain the inequality
$$
\frac{\d y}{\d t}\le\|f-g\|_m+\alpha y^2\qqwith
y(0)=\|u_0-v_0\|_m,
$$
where
\be{alphadef}
\alpha=c_m\exp\[c_m\int_0^T\(\|u(s)\|_m+\|u(s)\|_{m+1}\)\,\d s\].
\ee

In this case the condition (\ref{condition}) from lemma \ref{lemma}
becomes (\ref{condition2}). If this is satisfied then $y(t)$ is
uniformly bounded on $[0,T^*)$; the solution can therefore be
extended as a strong solution beyond $t=T^*$, contradicting the fact
that $T^*\le T$. It follows that $v(t)$ is a strong solution on
$[0,T]$ with the same regularity as $u$.
\end{proof}

\subsection{The Euler equations ($\nu=0$)}

We note here that the dissipative term $\nu Au$ plays no direct
r\^ole in the proof of theorem \ref{robust}; however, it does enter
indirectly via the regularity results of theorem \ref{reggie} that
are required to justify the equality in (\ref{langlerangle}). It is
possible to circumvent this via an appropriate mollification
(cf.~Constantin, E, \& Titi, 1994) and obtain the above result (and
those that follow) for solutions of the Euler equations ($\nu=0$) as
well as those of the Navier-Stokes equations.

\begin{theorem}\label{robust_euler}
Let $m\ge 3$ and let $u\in L^2(0,T;V^{m+3})$ be a solution of the 3d
Euler equations
$$
\frac{\d u}{\d t}+B(u,u)=f(t)\qqwith u(0)=u_0\in V^{m+3}
$$
where $f\in L^2(0,T;V^{m-1})\cap L^1(0,T;V^m)$. If $g\in
L^2(0,T;V^{m-1})$ and
\begin{eqnarray}
&&\|u_0-v_0\|_m+\int_0^T\|f(s)-g(s)\|_m\,\d
s\nonumber\\
&&\qquad\qquad<\frac{1}{c_mT}\exp\[-c_m\int_0^T
\(\|u(s)\|_m+\|u(s)\|_{m+1}\)\,\d s\]\label{condition2again}
\end{eqnarray}
then the equation
$$
\frac{\d v}{\d t}+B(v,v)=g(t)\qqwith v(0)=v_0\in V^{m+3},
$$
has a solution $v\in L^\infty(0,T;V^m)$.
\end{theorem}

\begin{proof}
Suppose that $\|v(t)\|_m<\infty$ for all $t\in[0,T^*)$ but that
$\|v(t)\|_m\rightarrow\infty$ as $t\rightarrow T^*$. As in the proof
of Theorem \ref{robust} we assume that $T^*\le T$ and obtain a
contradiction. First, observe that taking $T'<T^*$ it follows from
the regularity theory of Beale, Kato, \& Majda (1984) that $v\in
L^2(0,T';V^{m+3})$.

Now denote by $\varphi\in C^\infty_0(\Re^3)$ a standard mollifier
supported in the unit ball, set
$\varphi^\epsilon(x)=\epsilon^{-3}\varphi(x/\epsilon)$, and write
$u^\epsilon$ for the mollification of $u$ by convolution with
$\varphi^\epsilon$:
$$
u^\epsilon(x)=\int_{\Re^3}\varphi^\epsilon(y)u(x-y)\,\d y.
$$
(Note that the convolution is taken over all of $\Re^3$ with $u$ extended periodically, but the norms and inner products are still taken over $Q$.)

We return to (\ref{star}) with $\nu$ set to zero,
$$
\frac{\d w}{\d t}+B(u,w)+B(w,u)+B(w,w)=f-g,
$$
mollify the equation,
$$
\frac{\d w^\epsilon}{\d
t}+B(u,w)^\epsilon+B(w,u)^\epsilon+B(w,w)^\epsilon=(f-g)^\epsilon,
$$
and take the inner product with $w^\epsilon$ in $V^m$. In this way
we obtain (\ref{tooyoung}) with all quantities replaced by their
mollified counterparts, but with the addition of three error terms,
$$
\delta_{uww}+\delta_{wuw}+\delta_{www},
$$
where
$$
\delta_{uvw}=(B(u,v)^\epsilon,w^\epsilon)_m-(B(u^\epsilon,v^\epsilon),w^\epsilon)_m.
$$

Now, note that we have the pointwise identity (cf.~Constantin, E, \&
Titi, 1994)
$$
[(u\cdot\grad)v]^\epsilon=(u^\epsilon\cdot\grad)v^\epsilon+r_\epsilon(u,v)-[((u-u^\epsilon)\cdot\grad)(v-v^\epsilon)]
$$
where
$$
r_\epsilon(u,v)(x)=\int\varphi^\epsilon(y)[(u_y(x)\cdot\grad)v_y(x)]\,\d
y\qquad{\rm with}\qquad u_y(x)=u(x-y)-u(x).
$$
It follows that
$$
\|[(u\cdot\grad)v]^\epsilon-(u^\epsilon\cdot\grad)
v^\epsilon\|_m\le\|r_\epsilon(u,v)\|_m+\|((u-u^\epsilon)\cdot\grad)(v-v^\epsilon)\|_m.
$$
Since in 3d $|u(x)-u(y)|\le C\|u\|_{H^2}|x-y|^{1/2}$, it follows
that
$$
\|u-u^\epsilon\|_{H^m}\le
c\epsilon^{1/2}\|u\|_{H^{m+2}}\qquad\mbox{and}\qquad\|r_\epsilon(u,v)\|_{H^m}\le\Big(c^2\|u\|_{H^{m+2}}\|v\|_{H^{m+3}}\Big)\,\epsilon,
$$
and so we have
$$
|\delta_{uvw}|\le\Big(k\|u^\epsilon\|_{H^{m+2}}\|v^\epsilon\|_{H^{m+3}}\|w^\epsilon\|_{H^m}\Big)\epsilon.
$$

Now, use the fact that $\|u^\epsilon\|_{H^s}\le\|u\|_{H^s}$ for
$\epsilon<L/2$, we obtain
\bea
\ddt{}\|w^\epsilon\|_m&\le&
c_m\|u\|_m\|w^\epsilon\|_m+c_m\|u\|_{m+1}\|w^\epsilon\|_m+c_m\|w^\epsilon\|_m^2+\|f-g\|_m\\
&&+k'\Big(\|u\|_{m+2}\|w^\epsilon\|_{m+3}+\|w^\epsilon\|_{m+2}\|u\|_{m+3}+\|w^\epsilon\|_{m+2}\|w^\epsilon\|_{m+3}\Big)\epsilon
\eea
 for each $0<\epsilon<L/2$.

Now define (cf.~(\ref{alphadef}))
$$
\alpha_t=c_m\exp\(c_m\int_0^t\|u(s)\|_m+\|u(s)\|_{m+1}\,\d s\)
$$
to obtain
$$
\|w^\epsilon(t)\|_m\le
\frac{\eta_{T',\epsilon}}{1-\alpha_{T'}\eta_{T',\epsilon}T'}\qqfa
t\in[0,T']
$$
provided that
\bea
\eta_{T',\epsilon}&=&\|u_0-v_0\|_m+\int_0^{T'}\|f(s)-g(s)\|_m\,\d
s\\
&&\quad+\epsilon\,k'\int_0^{T'}\Big(\|u\|_{m+2}\|w^\epsilon\|_{m+3}+\|w^\epsilon\|_{m+2}\|u\|_{m+3}+\|w^\epsilon\|_{m+2}\|w^\epsilon\|_{m+3}\Big)\,\d
s\\
&<&\frac{1}{c_mT'}\exp\[-c_m\int_0^{T'}\(\|u(s)\|_m+\|u(s)\|_{m+1}\)\,\d
s\].
\eea
(cf.~(\ref{condition2})). Since both $u$ and $v$ are regular
solutions on $[0,T']$, it follows that
$$
\int_0^{T'}\|u\|_{m+2}\|w\|_{m+3}+\|w\|_{m+2}\|u\|_{m+3}+\|w\|_{m+2}\|w\|_{m+3}\,\d
s
$$
is finite. Given that $\|w^\epsilon\|_s\le\|w\|_s$ for all $\epsilon<L/2$ we can therefore let $\epsilon\rightarrow0$ and obtain
the bound
\be{wTp}
\|w^\epsilon(t)\|_m\le
\frac{\eta_{T'}}{1-\alpha_{T'}\eta_{T'}T'}\qqfa t\in[0,T']
\ee
where now
$$
\eta_{T'}=\|u_0-v_0\|_m+\int_0^{T'}\|f(s)-g(s)\|_m,
$$
provided that
\be{strictineq}
\eta_{T'}<\frac{1}{c_mT'}
\exp\[-c_m\int_0^{T'}\(\|u(s)\|_m+\|u(s)\|_{m+1}\)\,\d s\].
\ee

Since (\ref{condition2again}) holds, it follows that
(\ref{strictineq}) is verified for all $T'$ sufficiently close to
$T$. We can therefore deduce from (\ref{wTp}) that
$$
\|w(t)\|_m\le \frac{\eta_{T}}{1-\alpha_{T}\eta_{T}T}\qqfa t\in[0,T),
$$
and we have obtained a contradiction. It follows that $v$ is a
regular solution on $[0,T]$.\end{proof}

Similar techniques should be applicable to extend the other results
of this paper to the Euler case, but from now on we treat only the
Navier--Stokes equations.

\section{Deducing the existence of a strong solution via numerics}

As an application of theorem \ref{robust} -- and this may be the
most significant result in this paper -- we give an {\it a
posteriori} test to determine whether or not a numerical solution
of the 3d Navier-Stokes equations is meaningful. That is, we give
a criterion {\it depending only on the numerically computed
solution} that, if satisfied, guarantees that the exact equation
being approximated possesses a strong (and hence unique) solution.
Note that in the statement of the theorem the provenance of the
function $u$ is irrelevant.

\begin{corollary}\label{theretheyare}
Let $f\in L^2(0,T;V^m)$ and $u\in C^0([0,T];V^m)\cap
L^2(0,T;V^{m+1})$ with
$$
\frac{\d u}{\d t}+\nu Au+B(u,u)\ \in\ L^1(0,T;V^m)\cap
L^2(0,T;V^{m-1}).
$$
for some $m\ge3$. If $v_0\in V^m$ and
\begin{eqnarray}
&&\|u(0)-v_0\|_m+\int_0^T\left\|\frac{\d u}{\d t}(s)+\nu
Au(s)+B(u(s),u(s))-f(s)\right\|_m\,\d s\nonumber\\
&&\qquad\qquad<\frac{1}{c_mT}\exp\[-c_m\int_0^T\(\|u(s)\|_m+\|u(s)\|_{m+1}\)\,\d
s\]\label{gal2real}
\end{eqnarray}
then the solution of the Navier-Stokes equation
\be{exact}
\frac{\d v}{\d t}+\nu Av+B(v,v)=f(t)\qqwith v(0)=v_0\in V^m
\ee
is a strong solution on $[0,T]$ with $v\in L^\infty(0,T;V^m)\cap
L^2(0,T;V^{m+1})$.
\end{corollary}

\begin{proof}
The function $u$ is regular enough that it is the (unique) strong
solution of the Navier-Stokes equation (for the unknown $\tilde u$)
$$
\frac{\d \tilde u}{\d t}+\nu A\tilde u+B(\tilde u,\tilde
u)=\underbrace{\frac{\d u}{\d t}+\nu Au+B(u,u)}_{\mbox{notional
forcing}}\quad{\rm with}\quad \tilde u(0)=u_0;
$$
note that the conditions on $u$ ensure that the right-hand side is
an element of $L^2(0,T;V^{m-1})$. We now use theorem \ref{robust}
to compare $\tilde u$ with the solution of (\ref{exact}): in this
case the condition (\ref{condition2}) to guarantee that $v$ is a
strong solution is precisely (\ref{gal2real}).
\end{proof}

If $u$ comes from a discrete time-stepping algorithm, so that its
approximate values $u_n$ are only specified at times $t_n$ (with
$t_{n+1}>t_n$) then we can define a continuous function $u$ via
linear interpolation,
\be{interpolant}
u(t)=\frac{t_{n+1}-t}{t_{n+1}-t_n}\,u_n+\frac{t-t_n}{t_{n+1}-t_n}\,u_{n+1}\qquad\mbox{for}\qquad
t\in(t_n,t_{n+1}].
\ee
If $u_n\in V^{m+2}$ then certainly $u\in C^0([0,T];V^{m+2})$ and
$\d u/\d t\in L^2(0,T;V^{m+2})$. One can now apply the test of
corollary \ref{theretheyare} to this function $u$; if the test is
satisfied this again proves the existence of a strong solution for
(\ref{exact}).

\section{Convergence of the Galerkin approximations}

We now turn our attention to one particular form of numerical
solution, namely the Galerkin approximation. We show that given the
existence of a suitably smooth strong solution, this numerical
method is meaningful in that the Galerkin approximations can be
guaranteed to converge to the strong solution. Similar results --
convergence given the assumption that a strong solution exists --
are given for finite element methods by Heywood \& Rannacher (1982),
for a Fourier collocation method by E (1993), and for a nonlinear Galerkin method by Devulder, Marion, \& Titi (1993).

In some sense the result of Heywood (1982) that the Galerkin
approximations of a stable solution of the equations converge
uniformly on the whole time interval is in a similar spirit:
properties of the Galerkin method are deduced from an assumption on
the full equation. Of course, in this context it is perhaps more
natural to seek conditions under which one can guarantee the
existence of such a solution given properties of the Galerkin
approximations, see for example Constantin, Foias, \& Temam (1984)
and Titi (1987). For similar results for time-periodic solutions see
Titi (1991).

We should emphasise again that, in contrast to some related analyses of
the Galerkin method (e.g.~Rautmann, 1980) no assumption is made on the regularity of the
Galerkin approximations themselves.

\subsection{The Galerkin approximation}

We denote by $P_n$ the orthogonal projection in $H$ onto the first
$n$ eigenfunctions of the Stokes operator, and by $Q_n$ its
orthogonal complement. Denoting these eigenfunctions by
$\{w_j\}_{j=1}^\infty$, and their corresponding eigenvalues by $0<\lambda_1\le\lambda_2\le\cdots$, we have
$$
P_nu=\sum_{j=1}^n (u,w_j)w_j\qqand Q_nu=\sum_{j=n+1}^\infty
(u,w_j)w_j.
$$
Note that if $u\in V^m$ then
$$
\|Q_nu\|_m^2=\sum_{j=n+1}^\infty\lambda_j^m|(u,w_j)|^2\le\sum_{j=1}^\infty\lambda_j^m|(u,w_j)|^2=\|u\|_m^2
$$
and clearly $Q_nu\rightarrow0$ in $V^m$ as $n\rightarrow\infty$.

The Galerkin approximation of (\ref{nse_one}) is obtained by
projecting all terms onto the space $P_nH$:
\be{galerkin}
\frac{\d u_n}{\d t}+\nu Au_n+P_nB(u_n,u_n)=P_nf(t)\qqwith
u_n(0)=P_nu_0.
\ee

\subsection{Convergence of the Galerkin approximation}

Again we present our result for sufficiently strong solutions.
With some care one can combine the approach of Devulder et
al.~(1993) with that used here to give a proof for strong
solutions with minimal regularity; this will be presented
elsewhere.

\begin{theorem}\label{galerkin_thm} Let $u_0\in V^m$ with $m\ge 3$, $f\in
L^2(0,T;V^m)$, and let $u(t)$ be a strong solution of the
Navier-Stokes equations
\be{nse}
\frac{\d u}{\d t}+\nu Au+B(u,u)=f(t)\qqwith u(0)=u_0.
\ee
Denote by $u_n$ the solution of the Galerkin approximation
(\ref{galerkin}). Then $u_n\rightarrow u$ strongly in both
$L^\infty(0,T;V^m)$ and $L^2(0,T;V^{m+1})$ as
$n\rightarrow\infty$.
\end{theorem}

\begin{proof}
 The key, as with the robustness theorem,
is to arrange $B(u,u)-B(u_n,u_n)$ so that it only involves $u_n$
in the form $u_n-u$. Writing $w_n=u-u_n$ yields the equation
(cf.~(\ref{star}))
$$
\frac{\d w_n}{\d t}+\nu Aw_n+P_nB(u,w_n)+P_nB(w_n,u)+P_nB(w_n,w_n)=Q_nf-Q_nB(u,u).
$$
Taking the inner product of this equation with $A^mw_n$ we obtain
\begin{eqnarray}
\half\frac{\d}{\d
t}\|w_n\|_m^2+\nu\|w_n\|_{m+1}^2&\le&c_m\|u\|_m\|w_n\|^2+c_m\|u\|_{m+1}\|w_n\|_m^2+c_m\|w_n\|_m^3\nonumber\\
&&\qquad+\|Q_nf-Q_nB(u,u)\|_m\|w_n\|_m\label{forlater}.
\end{eqnarray}

Dropping the term $\nu\|w_n\|_{m+1}^2$ and dividing by $\|w_n\|_m$
yields
$$
\frac{\d}{\d t}\|w_n\|_m\le
c_m(\|u\|_m+\|u\|_{m+1})\|w_n\|+c_m\|w_n\|_m^2+\|Q_nf-Q_nB(u,u)\|_m.
$$
Setting
$$
y_n(t)=\|w_n(t)\|_m\exp\[-c_m\int_0^t\(\|u(s)\|_m+\|u(s)\|_{m+1}\)\,\d
s\]
$$
we obtain
$$
\dot y_n\le\|Q_nf-Q_nB(u,u)\|_m+\alpha y_n^2\qqwith
y_n(0)=\|Q_nu_0\|_m,
$$
where as in the proof of theorem \ref{robust}
$$
\alpha=c_m\exp\[c_m\int_0^T\(\|u(s)\|_m+\|u(s)\|_{m+1}\)\,\d s\].
$$

Noting that $\alpha$ is independent of $n$, and that $y_n(t)$ is
proportional to $\|w_n(t)\|_m$, convergence of the Galerkin
solutions will follow from convergence of $y_n(t)$ to zero. Using
lemma \ref{lemma} this will follow from
$$
\|Q_nu_0\|_m+\int_0^T\|Q_n[f(s)-B(u(s),u(s))]\|_m\,\d
s\rightarrow0\qas n\rightarrow\infty,
$$
which we now demonstrate.

That $\|Q_nu_0\|_m\rightarrow0$ as $n\rightarrow\infty$ is
immediate from the definition of $Q_n$. To show convergence of the
integral term, observe that since
\be{Buu}
\|B(u,u)\|_m\le c_m\|u\|_m\|u\|_{m+1}
\ee
and the regularity result of theorem \ref{reggie} guarantees that
$u\in L^2(0,T;V^{m+1})$, it follows that
$$
f(s)-B(u(s),u(s))\in V^m\qquad{\rm for}\qquad {\rm
a.e.}~s\in[0,T].
$$
We therefore know that $\|Q_n[f(s)-B(u(s),u(s))]\|_m$ converges
pointwise to zero for a.e.~$s\in[0,T]$, while it is clear that
$$
\|Q_n[f(s)-B(u(s),u(s)]\|_m\le\|f(s)-B(u(s),u(s))\|_m\quad{\rm
for}\quad {\rm a.e.}~s\in[0,T],
$$
and the right-hand side is an element of $L^1(0,T)$. It follows
from the Lebesgue dominated convergence theorem that
\be{integral}
\int_0^T\|Q_n[f(s)-B(u(s),u(s))]\|_m\,\d s\rightarrow0\qas
n\rightarrow\infty.
\ee

Returning to (\ref{forlater}) and integrating between $0$ and $T$
one can easily show the convergence of $w_n$ to zero in
$L^2(0,T;V^{m+1})$.
\end{proof}

It is clear that if the Galerkin approximations converge then one
can devise a fully discrete method that converges in a similar
sense, and this can easily be made precise. In the following
theorem one can take $u_{n,\delta}$ to be the linear interpolant
of a discrete set of values $u_{n,\delta}(t_j)\in V^m$ as in
(\ref{interpolant}) .

\begin{corollary}
Let $u_0\in V^m$ with $m\ge 3$, $f\in L^2(0,T;V^m)$, and let
$u(t)$ be a strong solution of the Navier-Stokes equations
(\ref{nse}).
For each $n$, let $\{u_{n,\delta}\}_{\delta>0}$ be a collection of
functions from $[0,T]$ into $P_nH$ with the property that
\be{finiteD}
u_{n,\delta}(t)\rightarrow u_n(t)\quad\mbox{in}\quad
L^\infty(0,T;V^m)\quad{\mbox as}\quad\delta\rightarrow0,
\ee
where $u_n$ is the solution of the Galerkin approximation
(\ref{galerkin}). Then there exist $\delta_n>0$ such that
$$
\sup_{0<\delta<\delta_n}\|u_{n,\delta}-u\|_{L^\infty(0,T;V^m)}
+\|u_{n,\delta}-u\|_{L^2(0,T;V^{m+1})}\rightarrow0
$$
as $n\rightarrow\infty$.
\end{corollary}

\begin{proof}
Since all norms on any finite-dimensional space are equivalent (in
particular those on $P_nH=P_nV^m=P_nV^{m+1}$) the convergence in
(\ref{finiteD}) also implies convergence in
$L^\infty(0,T;V^{m+1})$, and so in $L^2(0,T;V^{m+1})$. Otherwise
the result is immediate. \end{proof}

\section{Guaranteed numerical verification of the existence of a strong solution}

Our final result combines corollary \ref{theretheyare} and theorem
\ref{galerkin_thm} to show that the existence of a smooth strong
solution can be verified algorithmically {\it in a finite time} by
a sufficiently refined numerical computation. Namely, we show that
for $n$ sufficiently large the Galerkin solution {\it will}
satisfy the regularity test (\ref{gal2real}) from corollary
\ref{theretheyare}.

\begin{theorem}\label{verify}
Suppose that for some $m\ge 3$, $u_0\in V^m$, $f\in L^1(0,T;V^m)\cap
L^2(0,T;V^{m-1})$, and that $u$ is a strong solution of the
Navier-Stokes equations (\ref{nse}). Then there exists an $N$ such
that the solution $u_n$ of the Galerkin approximation
(\ref{galerkin}) satisfies condition (\ref{gal2real}) for every
$n\ge N$, i.e.~will pass the \emph{a posteriori} test for the
existence of a strong solution with data $(u_0,f)$.
\end{theorem}

\begin{proof}
First, we note that the convergence of the Galerkin approximations
to $u$ that is provided by theorem \ref{galerkin_thm} shows that
$\int_0^T\(\|u_n(s)\|_m+\|u_n(s)\|_{m+1}\)\,\d s$ is bounded
independently of $n$.

We need, therefore, only show that
$$
\|Q_nu_0\|_m+\int_0^T\left\|\frac{\d u_n}{\d t}(s)+\nu
Au_n(s)+B(u_n(s),u_n(s))-f(s)\right\|_m\,\d s
$$
(the left-hand side of (\ref{gal2real})) tends to zero as
$n\rightarrow\infty$. The requirement on the initial condition is
trivially satisfied, so we consider here only the `remainder term'
$$
\int_0^T\left\|\frac{\d u_n}{\d t}(s)+\nu
Au_n(s)+B(u_n(s),u_n(s))-f(s)\right\|_m\,\d s.
$$

Now, since $u_n$ satisfies the Galerkin approximation
(\ref{galerkin}) we have
$$
\frac{\d u_n}{\d t}(s)+\nu
Au_n(s)+B(u_n(s),u_n(s))-f(s)=Q_n[B(u_n(s),u_n(s))-f(s)],
$$
and so the remainder term is in fact
equal to
\be{newone}
\int_0^T \|Q_n[B(u_n(s),u_n(s))-f(s)]\|_m\,\d s.
\ee
In particular this shows that solutions of the Galerkin
approximation have the regularity required to ensure that the
integral on the left-hand side of (\ref{gal2real}) is
well-defined.

The integral in (\ref{newone}) coincides with the expression in
(\ref{integral}), whose convergence to zero we showed above in the
proof of theorem \ref{galerkin_thm}, except that the argument of
$B$ is $u_n$ rather than $u$. However,
$$
B(u_n(s),u_n(s))-B(u(s),u(s))=B(u_n(s)-u(s),u_n(s))+B(u(s),u_n(s)-u(s))
$$
and so
\bea
&&\int_0^T\|Q_n[B(u_n(s),u_n(s))-B(u(s),u(s))\|_m\\
&&\qquad\le\int_0^T\|B(u_n(s),u_n(s))-B(u(s),u(s))\|_m\,\d s\\
&&\qquad\le\int_0^T\|B(u_n(s)-u(s),u_n(s))\|_m\,\d s
+\int_0^T\|B(u(s),u_n(s)-u(s))\|_m\,\d s\\
&&\qquad\le c_m\int_0^T\|u_n(s)-u(s)\|_m\|u_n(s)\|_{m+1}\,\d
s\\
&&\qquad\qquad\qquad\qquad+c_m\int_0^T\|u(s)\|_m\|u_n(s)-u(s)\|_{m+1}\,\d
s
\eea
Since $u_n\rightarrow u$ strongly in both $L^\infty(0,T;V^m)$ and
$L^2(0,T;V^{m+1})$ the result follows.
\end{proof}

Once more it is possible to treat fully discrete schemes within a
similar framework. Supposing as above that a scheme gives rise to a
discrete set of values $u_{n,\delta}(t_j)$, one can define a linear
interpolation $u_{n,\delta}$. We then have the following result:

\begin{corollary}
Suppose that for some $m\ge 3$, $u_0\in V^m$, $f\in L^1(0,T;V^m)\cap
L^2(0,T;V^{m-1})$, and that $u$ is a strong solution of the
Navier-Stokes equations (\ref{nse}). For each $n$, let
$\{u_{n,\delta}\}_{\delta>0}$ be a collection of functions from
$[0,T]$ into $P_nH$ with the property that
\be{finiteDt}
\d u_{n,\delta}(t)/\d t\rightarrow \d u_n(t)/\d t\quad\mbox{in}\quad
L^\infty(0,T;V^m)\quad{\mbox as}\quad\delta\rightarrow0,
\ee
where $u_n$ is the solution of the Galerkin approximation
(\ref{galerkin}). Then there exists an $N$ and a sequence $\delta_n$
such that the interpolant $u_{n,\delta}$ satisfies condition
(\ref{gal2real}) for every $n\ge N$ and $\delta<\delta_n$, i.e.~the
fully discrete numerical solution will pass the \emph{a posteriori}
test for the existence of a strong solution with data $(u_0,f)$.
\end{corollary}

\begin{proof}
The interpolant $u_{n,\delta}$ satisfies the Galerkin
approximation (\ref{galerkin}) except for an error
$$
r_{n,\delta}=\frac{\d u_n}{\d t}-\frac{\d u_{n,\delta}}{\d t};
$$
this error converges to zero in $L^\infty(0,T;V^m)$ because of
(\ref{finiteDt}) (cf.~Higham \& Stuart, 1998). It follows that we
now obtain
$$
\int_0^T \|Q_n[B(u_n(s),u_n(s))-f(s)+r_{n,\delta}(s)]\|_m\,\d s.
$$
rather than (\ref{newone}), where $r_{n,\delta}$ converges to zero
in $L^\infty(0,T;V^m)$. This is clearly sufficient to follow the
argument in the proof of theorem \ref{verify}.
\end{proof}

\section{Conclusion}

Despite the lack of a guarantee that unique solutions exist for the
three-dimensional Navier-Stokes equations, we have shown that it is
possible to perform `rigorous' numerical experiments. In particular,
we have given an {\it a posteriori} test that, if satisfied by a
numerical solution, guarantees that it approximates a true strong
solution of the Navier-Stokes equations. Remarkably, the existence
of such a solution can be verified using such numerical computations
in a finite time; some computations along these lines will be reported in a future publication.

\section*{Acknowledgments}
This work arose from discussions at the workshop on {\it
Singularities, coherent structures and their role in intermittent
turbulence} at Warwick Mathematics Institute in September 2005,
which was supported by the EPSRC. SIC's travel and workshop fee were
funded from EPSRC grants GR/S67029 and GR/S82947. The work of PC was
partially supported by the National Science Foundation, grant no.
DMS-0504213. JCR is a Royal Society University Research Fellow and
would like to thank the society for all their support. The work of
EST was supported in part by the National Science Foundation grant
no. DMS-0504619, the BSF grant no. 200423, the USA Department of
Energy under contract W-7405-ENG-35, and the ASCR Program in Applied
Mathematical Sciences. We would also like to thank the referees for
their helpful comments.

\bigskip

\parindent0pt

\begin{small}

Constantin, P. 1986 Note on loss of regularity for solutions of the
$3$-D incompressible Euler and related equations. {\it Comm. Math.
Phys.} {\bf 104}, 311--326

Constantin, P., E, W., \& Titi, E.S. 1994 Onsager's conjecture on
the energy conservation for solutions of Euler's equation. {\it
Commun. Math. Phys.} {\bf 165}, 207--209.

Constantin, P. \& Foias, C. 1988 {\it Navier-Stokes equations}.
Chicago: University of Chicago Press.

Constantin, P., Foias, C., \& Temam, R. 1984 On the large time
Galerkin approximation of the Navier-Stokes equations. {\it SIAM J.
Numer. Anal.} {\bf 21}, 615--634.

Devulder, C., Marion, M., \& Titi, E.S. 1993 On the rate of convergence of the nonlinear Galerkin methods. {\it Math. Comp.} {\bf 60}, 495--514.

E, W. 1993 Convergence of Fourier methods for the Navier--Stokes
equations. {\it SIAM J. Num. Anal.} {\bf 30}, 650--674.

Foias, C., \& Temam, R. 1989 Gevrey class regularity for the solutions of the Navier--Stokes equations. {\it J. Func. Anal.} {\bf 87}, 359--369.

Heywood, J.G. 1982 An error estimate uniform in time for spectral
Galerkin approximations of the Navier-Stokes problem. {\it Pacific
J. Math.} {\bf 98}, 333--345

Heywood, J.G. \& Rannacher, R. 1982 Finite element approximation of
the nonstationary Navier-Stokes problem. I. Regularity of solutions
and second-order error estimates for spatial discretization. {\it
SIAM J. Numer. Anal.} {\bf 19}, 275--311.

Higham, D.J. \& Stuart, A. M. 1998 Analysis of the dynamics of local
error control via a piecewise continuous residual. {\it BIT} {\bf
38}, 44--57.

Ponce, G., Racke, R., Sideris, T.C., \& Titi, E.S. 1993 Global
stability of large solution to the $3-D$ Navier-Stokes equations.
{\it Comm. Math. Phys.} { \bf 159}, pp.~329--341.

Rautmann, R. 1980 On the convergence rate of nonstationary
Navier-Stokes approximations. In {\it Approximation methods for
Navier-Stokes problems} (ed. R. Rautmann). Springer Lecture Notes in
Mathematics, no. 771, pp.~425--449.

Robinson, J.C. 2001 {\it Infinite-dimensional dynamical systems}.
Cambridge: Cambridge University Press.

Temam, R. 1977 {\it Navier-Stokes Equations, Theory and Numerical
Analysis}. Amsterdam: North Holland. Revised 1984 version reprinted
by AMS Chelsea, 2001.

Titi, E.S. 1987 On a criterion for locating stable stationary
solutions to the Navier-Stokes equations. {\it Nonlinear Anal.} {\bf
11}, 1085--1102.

Titi, E.S. 1991 Un Crit\`ere pour l'approximation des solutions
p\'eriodiques des \'equations de Navier--Stokes. {\it Comptres
Rendus de L'Acad\'emie des Sciences, Paris, S\'erie I} {\bf 312},
pp.~41--43.

\end{small}

\end{document}